\documentclass[12pt,a4paper,twoside,leqno]{article}
\usepackage{amsthm,amsmath,amssymb,mathabx}
\usepackage{color}
\newtheorem{teor}{Theorem}[section]
\newtheorem{defin}[teor]{Definition}
\newtheorem{lemm}[teor]{Lemma}
\newtheorem{osse}[teor]{Remark}
\newtheorem{prop}[teor]{Proposition}
\newtheorem{defi}[teor]{Definition}
\newtheorem{coro}[teor]{Corollary}
\newtheorem{prob}[teor]{Problem}
\newtheorem{hypo}[teor]{Hypothesis}

\newcommand{\bele}{\begin{lemm}\begin{sl}}
\newcommand{\enle}{\end{sl}\end{lemm}}
\newcommand{\bedef}{\begin{defi}\begin{sl}}
\newcommand{\eddef}{\end{sl}\end{defi}}
\newcommand{\bete}{\begin{teor}\begin{sl}}
\newcommand{\ente}{\end{sl}\end{teor}}
\newcommand{\beos}{\begin{osse}\begin{sl}}
\newcommand{\eddos}{\end{sl}\end{osse}}
\newcommand{\bepr}{\begin{prop}\begin{sl}}
\newcommand{\empr}{\end{sl}\end{prop}}
\newcommand{\bepro}{\begin{prob}\begin{rm}}
\newcommand{\empro}{\end{rm}\end{prob}}
\newcommand{\bede}{\begin{defin}\begin{sl}}
\newcommand{\edde}{\end{sl}\end{defin}}
\newcommand{\beco}{\begin{coro}\begin{sl}}
\newcommand{\enco}{\end{sl}\end{coro}}
\newcommand{\behy}{\begin{hypo}\begin{sl}}
\newcommand{\enhy}{\end{sl}\end{hypo}}

\newcommand{\prova}{{\bf Proof.\hspace{4mm}}}

\newcommand{\thspace}{\hspace{3mm}}

\newcommand{\RR}{\mathbb{R}}
\newcommand{\NN}{\mathbb{N}}

\newcommand{\beeq}[1]{\begin{equation}\label{#1}}
\newcommand{\eddeq}{\end{equation}}

\newcommand{\beeqa}[1]{\begin{eqnarray}\label{#1}}
\newcommand{\eddeqa}{\end{eqnarray}}

\newcommand{\beal}[1]{\begin{align}\label{#1}}
\newcommand{\eddal}{\end{align}}

\newcommand{\bespl}[1]{\begin{split}\label{#1}}
\newcommand{\edspl}{\end{split}}

\newcommand{\bega}[1]{\begin{gather}\label{#1}}
\newcommand{\edga}{\end{gather}}

\newcommand{\beeqax}{\begin{eqnarray*}}
\newcommand{\eddeqax}{\end{eqnarray*}}

\def\qed{\ifmmode % if math mode, assume display: omit penalty etc.
  \else \leavevmode\unskip\penalty9999 \hbox{}\nobreak\hfill
  \fi
  \quad\hbox{\hskip.5em\vrule width.4em height.6em depth.05em\hskip.1em}}
\def\endproofsym{\qed}

\def\endnobox{\def\endproofsym{}\end{proof}\def\endproofsym{\qed}}

\newcommand{\no}{\nonumber}

\newcommand{\beeqao}{\begin{eqnarray}\no}
\newcommand{\bealo}{\begin{align}\no}
\newcommand{\besplo}{\begin{split}\no}
\newcommand{\begao}{\begin{gather}\no}

\newcommand{\duav}[1]{\langle{#1}\rangle}

\newcommand{\dt}{\partial_t}

\newcommand{\io}{\int_\Omega}

\newcommand{\e}{\varepsilon}

\newcommand{\dd}{\, \mathrm{d}}

\DeclareMathOperator{\dive}{div}

\newcommand{\weak}{\mathop\rightharpoonup}
\newcommand{\weakstar}{\mathop{\rightharpoonup}^{*}}

\newcommand{\chin}{\chi_\e}

\newcommand{\sn}{s_\e}

\newcommand{\BBB}{\color{black}}  
\newcommand{\BBBnew}{\color{black}}  
\newcommand{\EEE}{\color{black}}
\newcommand{\EEEbis}{\color{black}}

\newcommand{\EEEtris}{\color{black}}

\newcommand{\BL}{\color{black}}

\numberwithin{equation}{section}
\begin{document}

\title{Generalized gradient flow structure \\ of \EEEtris internal \BBB energy \BL driven phase field systems}

%\date{}

\author{Elena Bonetti\\
{\sl Dipartimento di Matematica, Universit\`a degli Studi di Pavia}\\
{\sl Via Ferrata, 1, I-27100 Pavia, Italy}\\
{\rm E-mail:~~\tt elena.bonetti@unipv.it}\\
\and
Elisabetta Rocca\\
{\sl Weierstrass Institute for Applied
Analysis and Stochastics}\\ 
{\sl Mohrenstr.~39, D-10117 Berlin, Germany}\\ 
{\sl and}\\
{\sl Dipartimento di Matematica, Universit\`a degli Studi di Milano}\\ 
{Via Saldini 50, I-20133 Milano, Italy}\\
{\rm E-mail:~~\tt elisabetta.rocca@wias-berlin.de}\\
 {\rm and \tt  elisabetta.rocca@unimi.it}}
\maketitle

\begin{abstract}
\BBB
In this paper we introduce a general abstract formulation of a variational thermomechanical model, 
by means of a unified derivation via \EEEtris a generalization of the \BBB principle of virtual powers for all the variables of the system, including the 
thermal one. In particular, choosing as thermal variable the entropy of the system, and as driving functional the internal energy,  
we get a gradient flow structure (in a suitable abstract setting) 
for the whole nonlinear PDE system. We prove a global in time existence of (weak) solutions result for the Cauchy problem 
associated to the abstract PDE system as well as uniqueness in case of suitable smoothness assumptions on the functionals.  
\BL
\end{abstract}

\vspace{.2cm}

\noindent
{\bf Key words:}\thspace \BBB gradient flow, phase field systems, existence of weak solutions, uniqueness. 
\BL

\vspace{2mm}

\noindent
{\bf AMS (MOS) subject clas\-si\-fi\-ca\-tion:}\thspace \BBB 74N25, 82B26, 35A01, 35A02. \BL

\vspace{2mm}

%\newpage

%%%%%%%%%%%%%%%%%%%%%%%%%%%%%%%%%%%%%%%%%%%%%%%%%%%%%%%%%%%%%%%%%%%
%%%%%%%%%%%%%%%%%%%%%% Caption da ambo i lati %%%%%%%%%%%%%%%%%%%%%

\pagestyle{myheadings}
\newcommand\testopari{\sc Bonetti -- Rocca}
\newcommand\testodispari{\sc \BBB Energy \BL driven systems}
\markboth{\testodispari}{\testopari}

%%%%%%%%%%%%%%%%%%%%%%%%%%%%%%%%%%%%%%%%%%%%%%%%%%%%%%%%%%%%%%%%%%%
%%%%%%%%%%%%%%%%%%%%%%%%%%%%%%%%%%%%%%%%%%%%%%%%%%%%%%%%%%%%%%%%%%%

%%%%%%%%%%%%%%%%%%%%%%%%%%%%%%%%%%%%%%%%%%%%%%%%%%%%%%%%%%%%%%%%%%%
%%%%%%%%%%%%%%%%%%%%%%%%% Introduction %%%%%%%%%%%%%%%%%%%%%%%%%%%%

\section{Introduction}
\label{intro}

In this paper we introduce a general derivation of thermo-mechanical phase transition models by use of a generalization of the principle of virtual powers, in which micro-forces and thermal forces \BBB are \BL  included.  It is known that a recent field of research, in the framework of phase transitions, has concerned models with \BBB some \BL micro-forces (see, e.g., the approaches by Fr\'emond \BBB \cite{fremondlibro} \EEEbis and  \BBB by Gurtin \BL \cite{gurtin}\BL). \EEEbis The main idea is that the equations governing the evolution of phase transition phenomena may be derived by a variational principle, i.e. the principle of virtual powers, in which micro-forces, responsible for phase transitions (i.e. for changes in the microstructure level of the materials), are included. As a consequence, the resulting PDE system provides an intrinsic variational structure, at least concerning equations for displacements and internal quantities, as phase or order parameters. Many authors have dealt with this kind of approach. We mention, among the others \BBB we quote some contributions as \cite{lss}, \cite{fpr}, and \cite{ms}. \BL 

On the other hand, as far as thermal properties \BBB are concerned\BL, in the recent years several efforts have been spent to investigate models in which an entropy balance (or imbalance) equation was introduced in place of the more classical ``heat equation''. We recall, e.g., the contribution by \cite{bcf}, \cite{BFI}, and \cite{BFR}. In particular, let us mention that the last paper shows a derivation of the equation on the entropy by convex analysis tools and the application of a Legendre transformation for the free energy. It is interesting to observe that in this framework, also thermal memory is formally justified from the point of view of the derivation of the model.
\BL 
In a different direction Podio-Guidugli, in relation to a theory proposed by Green and Naghdi, introduced the possibility of including thermal displacements and forces in the whole balance of the principle of virtual powers, so that the entropy equation may be recovered, as well as the momentum equation, as a  ``balance of forces'', forcing the system on the base of some ``reluctance to order''. 
\EEEbis 
Indeed, starting from the consideration that some virtual power principle may be used to deduce all balance and imbalance laws of thermomechanics, he suggested to use it also for the derivation of thermal evolution, through the notion of  thermal displacement. As a consequence, he derives an equation for the entropy of the system, which is combined with momentum balance. This approach turns out to be consistent with thermodynamical principles. See, among the others, \cite{podio} and the papers by Green and Naghdi 
 \cite{green1}, \cite{green2}, and references therein. 
\BBB 
Finally, we can quote the recent contribution \cite{mielke}, where a gradient structure of systems in thermoplasticity is introduced by means of a free entropy functional instead of the internal energy, which is the driving functional in the present contribution. 
\BL 

Indeed, in this paper, we aim to 
\BL combine the previous approaches and provide a general abstract formulation of a  variational thermomechanical model which can be applied to recover different phase transitions and phase separation phenomena, also accounting for mechanical or thermal effects. \EEEbis Hence, we introduce an unified approach which formally justifies the evolution \BBB of the thermal variable (represented here by the entropy of the system), \BL  the phase parameters, and (possibly) the displacements. Actually, in the following, we are dealing with two state variables:  $s$, which mainly plays the role of the entropy, and a phase parameter $\chi$, \BBB representing the internal mechanical variable.  The main advantage of the gradient structure is the possibility of deriving  a time-incremental minimization procedure, where the internal energy functional is minimized with respect to the entropy and the internal variables and so the existence of weak solutions for the associated Cauchy problem can be deduced under quite general assumptions on the involved nonlinearities. \BL 

Indeed, the choice of the energy functional and the dissipation potential  are \BL fairly general. In particular, in the internal energy functional we can include  multivalued operators to ensure some internal constraints. \BBB Since the resulting \BBB gradient flow \BL structure is nonlinear \BBB and non smooth, \BL we have to introduce a suitable notion of (weak) solution \BBB in order to get  a global in time existence result. \BL  However, the weak notion of solution we are introducing is naturally in accordance with the physical meaning of the problem \BBB under consideration as well as with the classical principles of thermodynamics. 
The proof is performed by means of a combined regularization and time discretization procedure. 
Moreover, uniqueness of solutions is proved under some further smoothness assumption on the internal energy functional. \BL 

The paper is organized as follows. In the next Section~\ref{model} we derive the model and state the main assumptions on the involved physical quantities and functionals. The main existence result is stated and proved in Section~\ref{sec:main}, as well as the uniqueness of solutions. 
\BL

\section{The model and the main assumptions}
\label{model}

Let $\Omega\subset\RR^3$  be a bounded and sufficiently regular domain with boundary $\Gamma:=\partial\Omega$. We introduce an Hilbert triplet 
$V\subset H\subset V'$
(with dense and compact injections), where $V=H^1(\Omega)$, $H=L^2(\Omega)$, and $H$ is identified as usual  with its dual. We introduce the notations  \BL $\duav{\cdot, \cdot}$
for  the duality pairing between $V$ and $V'$  and $(\cdot,\cdot)$  for the usual
scalar product both in $H$ and in $L^2(\Omega)^3$. To simplify  the notation, we  write  $H$ in place
of $L^2(\Omega)^3$, or $V$ in place of $H^1(\Omega)^3$, when vector-valued functions are
considered. For every $f \in V'$ we
indicate by $\overline{f}$ the {\it spatial mean} of $f$ over $\Omega$, i.e.
\[
  \overline{f} := \frac{1}{|\Omega|} \duav{ f, 1},
\]
where  $|\Omega|$ stands for the Lebesgue measure of $\Omega$.
We note as $H_0$, $V_0$ and $V_0'$ the closed subspaces of functions
(or functionals) having zero mean value in $H$, $V$, and  in $V'$, respectively. 
Then, by the Poincar\'e-Wirtinger inequality,
\[
  \| v \|_{V_0} := \left( \io | \nabla v |^2 \,{\rm d}x \right )^{1/2}
\]
represents a norm on $V_0$ which is equivalent to the norm naturally  inherited from $V$. In particular $\|\cdot\|_{V_0}$ is a Hilbert norm  \EEEbis associated to a scalar product  $((\cdot,\cdot))_{V_0}$  (defined in \eqref{defiJ}), and thus \BL
we can introduce the associated
Riesz isomorphism mapping $A:V_0 \to V_0'$ by setting, for $u,v\in V_0$,
\begin{equation}\label{defiJ}
  \duav{A u, v}:= (\!(u,v)\!)_{V_0}
   :=\int_\Omega\nabla u\cdot \nabla v\, {\rm d} x,
\end{equation}
so that $\duav{Au, u}=\|u\|_{V_0}^2$ for every $u\in V_0$ and $\duav{v, A^{-1}(v)}=\|v\|_{V_0'}^2$ for every $v\in V_{0}'$. \BL 
Finally, we can identify $H_0$ with $H_0'$
by means of the scalar product of $H$ so to obtain the Hilbert
triplet $V_0 \subset H_0 \subset V_0'$, where inclusions are continuous
and dense. In particular, if $z\in V$ and $v \in V_0$,
it is easy to see that
\begin{equation}\label{duale}
  \io \nabla z \cdot \nabla (A^{-1} v) \, {\rm d}x
   = \io ( z - \bar z ) v \, {\rm d}x
   = \io z v \, {\rm d}x.
\end{equation}
In what follows in  this section we introduce our modelling approach and the set of PDEs and initial and boundary conditions which we are going to analyze in the next sections. \BL

\subsection{The Principle of Virtual Powers}  

The model is derived by using a variational principle in mechanics which is known as (generalized) principle of virtual power. Indeed, we refer to some generalization of the well known mechanical principle as we are including in the involved forces  the  microscopic  forces, acting on some ``micro-scale'', and also possible ``thermal forces". Without entering the details of this argumentation,  let us  point out that this principle is formally based on the fact that velocities are considered in a suitable linear space and thus forces are defined as elements acting on velocities with respect to some duality relation between the two spaces. This is done for any (sufficiently smooth) subdomain ${\cal D}\subseteq\Omega$. Hence, before proceeding we make precise the virtual velocities we are considering. More precisely, let us take the couple of virtual velocities $(\delta_t,v_t)$ (whose physical meaning may change time to time). In the case when no accelerations are included, the principle of virtual powers can be written  considering the power of internal forces $\mathcal{P}_{int}$ and of external forces $\mathcal{P}_{ext}$ \BBB (depending on ${\cal D}$, $\delta_t$, $v_t)$ as follows:\BL
\[
\mathcal{P}_{int}({\cal D},\delta_t, v_t)+\mathcal{P}_{ext}({\cal D},\delta_t, v_t)=0.
\]
We assume that  the power of internal forces is introduced as follows (in $\Omega$ and for any virtual velocities $\delta_t\in V_0$ and $v_t\in V$)
\begin{align}\label{servevariazionale}
\mathcal{P}_{int}&=\langle\langle F, \delta_t\rangle\rangle+\langle\langle{\cal G},v_t\rangle\rangle=\duav{F, \delta_t}+\io B\,v_t\dd x+\io E\,\nabla v_t\dd x,
\end{align}
%%%%TAGLIARE
where $F$, $B$ and $E$ denote interior thermal and mechanical (micro) forces and stresses, respectively and the duality relation $\langle\langle\cdot,\cdot\rangle\rangle$ is suitably defined between forces and velocities spaces. 

 Analogously, the power of external forces is
 $$
 \mathcal{P}_{ext}=\langle\langle{\cal R},\delta_t\rangle\rangle+\langle\langle{\cal Z} ,v_t\rangle\rangle.
 $$
We let $v_t\in V$ and $\delta_t\in V_0$ and we  assume there exists  $Z,$ $z$  such that
\begin{equation}\label{RZ}
\langle\langle{\cal R},\delta_t\rangle\rangle=\langle{\cal R},\delta_t\rangle\quad \hbox{and}\quad\duav{{\cal Z}, v_t\BL}=\io Z \,v_t+\int_\Gamma z \,v_t,
\end{equation}
where $Z$ and $z$ stand for  the external forces acting in the bulk $\Omega$ and at the boundary $\Gamma$, respectively.

It is clear that that we are considering a different behavior of the forces on the two types of virtual velocities. Indeed, we note that the elements ${\cal G}$ and ${\cal Z}$ are defined as a.e. forces  living  in the bulk and on the boundary (with suitable summability), while we take $F$  and ${\cal R}$ as general as possible to include all the different (and less regular) situations we will face. In particular, as it will be clear once we will make a precise choice of the actual velocities (cf.~\eqref{vps} and \eqref{vpchi1}), of the energy functional \eqref{inten}, and of the dissipation potential \eqref{phi}, we aim to write down an equation for the thermal variable of {\sl conservative} type: it will result indeed as a {\sl conservation of energy}, while the equation for the mechanical variable will be on {\sl non conservative type}. This mainly motivates the choice we have made for the power of internal and external forces. Other choices are possible (cf., e.g., Remark~\ref{case2}), but we prefer not to move in this direction in the present contribution.

\subsection{The constitutive relations and the PDEs}
\label{state-en-diss}

\paragraph{The state variables.} 
We are dealing with a physical system governed by the state variables $(s,\chi,\nabla \chi)$ whose evolution is ruled by different thermomechanical relations. 
Note that we are distinguishing between different {\sl dependence} of the energy with respect to the   two variables  $s$ and $\chi$: we consider, in particular, the gradient $\nabla\chi$ but not $\nabla s$ as state variable  (cf.~\eqref{inten}).
This corresponds to the specific choice we have done for the forces $F$, ${\cal G}$ and ${\cal R}, {\cal Z}$  we have made in \eqref{servevariazionale} and \eqref{RZ}. 

In order to get \BL the evolution of $s$, we take the actual velocities as $\delta_t=A^{-1}(\xi_t)$, where $\xi_t\in V_0'$ and $v_t=0$, in order  to get (for ${\cal R}=0$) $\duav{F, A^{-1}(\xi_t)}=0$ for all $\xi_t\in V_0'$ and so we obtain 
\begin{equation}\label{vps}
A(F)=0\quad\hbox{ in }V_0'.
\end{equation}
The evolution of $\chi$ is obtained by integrating by parts in $\mathcal{P}_{int}$ and choosing $Z=z=0$ and $\delta_t=0$ as well\BL:
\begin{equation}\label{vpchi1}
B-\dive E=0 \quad \hbox{in }V',\end{equation}\BL
with the   no-mass flux through the boundary of $\Omega$\BL :
\begin{equation}\label{vpbouchi}
E\cdot {\bf n}=0\quad\hbox{on }\Gamma\,,
\end{equation}
where we have denoted by ${\bf n}$ the outward unit normal vector to $\Gamma$. \BL 
Notice that we have chosen here to have 0 ezternal forces ${\cal R}$ and ${\cal Z}$  only for simplicity of notation.\EEEbis
\beos
In the following, we mainly refer to the variable $\chi$ as a phase or order parameter, i.e. \BBB related to the micro-structure of the physical system. \BL However, let us point out that we  could formally include in our procedure the derivation of the (more) classical momentum balance equation (letting, e.g., $\chi$ stand for displacements). In this case,  the force $B$ has to be equal to 0, due to the principle of rigid motions. \eddos

\paragraph{The functionals and the main assumptions.} We introduce two functionals governing the evolution and the equilibrium of our (thermo)mechanical system. These functionals depend on the state variables and on the dissipative variables, respectively. As far as the equilibrium, it is governed by an energy functional, and we choose to make use of an internal energy functional (in place of the free energy functional). This choice is motivated by the fact that we may interprete $s$ as the entropy of the system (see \cite{mielke} for a  physical justification). However, it is well known that, under suitable assumptions, the internal energy may be introduced as the Legendre transformed of the free energy. \BL

Before we make precise the choices of the  internal energy functional and of the dissipation potential, let us introduce  a  function $W$, \EEEbis depending on $\chi$,  as the sum of a convex possibly non-smooth part and non-convex but regular function and it  satisfies some smoothness and growth assumptions, in particular, we need:
\behy\label{HW} Assume $W(\chi)=\widehat\beta(\chi)+\widehat\gamma(\chi)$, where 
\begin{description}
\item[{\bf (w1)}] $\widehat\beta: {\rm dom}(\widehat\beta)\to [0,+\infty]$ is convex, proper, and lower semicontinuous, 
\item[{\bf (w2)}] $\widehat\gamma:\RR\to\RR$ is a $C^{1,1}$ function on $\RR$.
\end{description} 
\enhy

\beos\label{exW}
Notice that particularly meanigful choices of $W$  used in the literature of phase transitions (when  $\chi$ assumes the meaning of phase variable) \BL are,  for example, \BL
\begin{itemize}
\item[1.] the {\sl double well potential} $W(\chi)=(\chi^2-1)^2$ 
\item[2.] the {\sl logarithmic potential} $W(\chi)= \chi\log(\chi)+(1-\chi)\log(1-\chi)-\chi^2$
\item[3.] the {\sl double obstacle potential} $W(\chi)=I_{[0,1]}(\chi)-\chi^2$, where $I_{[0,1]}$ denotes the indicator function of the interval $[0,1]$  and it is defined as $I_{[0,1]}(x)=0$ if $x\in[0,1]$ and $I_{[0,1]}(x)=+\infty$ otherwise. \BL
\end{itemize} 
\eddos
Moreover, we introduce a function $j(\theta,\chi)\,:\,\RR\times \RR\to [0,+\infty]$ such that 
\begin{align}\no
&\theta \mapsto j(\theta, \chi) \quad \hbox{ is a convex, proper, and lower semicontinuous  for every $\chi\in \RR$ and}\\
\no
& \chi\mapsto j(\theta, \chi) \quad \hbox{ is a $C^1$ function for every }\theta\in \RR, 
\end{align}
and let 
\begin{align}\label{JH}
&J_H(\theta,\chi)=\begin{cases}\displaystyle\int_\Omega j(\theta,\chi) \quad &\hbox{if } (\theta,\chi)\in H\times H\quad\hbox{and }j(\theta,\chi)\in L^1(\Omega)\\
+\infty\quad &\hbox{if } (\theta,\chi)\in H\times H\quad\hbox{and }j(\theta,\chi)\varnotin L^1(\Omega)
\end{cases} \\
&
\no
\\
\label{JV}
&J_V(\theta,\chi)=J_H(\theta,\chi)\quad\hbox{\BBB on \BL } V\times H.
\end{align}
Hence, we can introduce the convex conjugate of $J_V$ as follows 
$J_V^*(s,\chi)\,:\, V'\times H\to [0,+\infty]$ is defined as 
\begin{equation}\label{JstarV}
J_V^*(s,\chi)=\sup_{\theta\in V}\left(\duav{s, \theta}-J_V(\theta,\chi)\right), \quad (s,\chi)\in V'\times H.
\end{equation} 
Now, we are in the position of introducing \BL the energy functional $e:V'\times H\times H\to (-\infty,+\infty]$\BL :
\begin{equation}\label{inten}
e(s,\chi,\nabla\chi)=J_V^*(s,\chi)+\io\left(\frac12|\nabla\chi|^2+W(\chi)\right)\dd x.
\end{equation}
Let us note here that the first term in \eqref{inten}  contains both the purely caloric part of the energy functional (i.e.~the one depending only on $s$ as well as the coupling terms depending on both $s$ and $\chi$) (cf.~Subsection~\ref{ex:jstar} for possible choices of $J_V^*$), while inside the integral over $\Omega$ we have the parts accounting for the {\sl nonlocal} interfacial energy effects (the  $|\nabla\chi|^2$) and the mixing potential $W$ (cf.~Remark~\ref{exW} for examples of possible choices of functions $W$). We intentionally choose not to consider interfacial ({\sl nonlocal}) energy effects in the variable $s$ in order to differentiate the roles of the caloric and the mechanical parts ($s$ and $\chi$, respectively) in our approach. 

Then, we can define the subdifferential (with respect to the variable $s$)  $\partial_{V',V} J_V(s,\chi)$ which maps $V'\times H$ into $2^V$ as \BBB (cf., e.g., \cite{brezis}):\BL 
\begin{align}\label{partialV}
&v\in \partial_{V',V}J_V^*(s,\chi) \quad \hbox{ if and only if } v\in V, \\
\no
&\qquad (s, \chi)\in D(J_V^*),\hbox{ and } J_V^*(s,\chi)\leq \duav{s-w,v}+J_V^*(w,\chi)\quad \forall (w,\chi)\in V'\times H\,.
\end{align}

Actually, in what follows we will always work in the space $V'\times H$ and so we will state directly the assumptions we need on the functional $J_V^*$ defined in \eqref{JstarV}. In particular, we need the following assumptions:
\behy~\label{HJ} We assume that \BBB $J_V^*:V'\times H\to [0,+\infty]$ is such that: \BL there exist two positive constants $c_1,c_2\in \RR^+$ such that the functional $J_V^*$ defined in \eqref{JstarV} satisfies:
\begin{description}
\item[(J1)] $\chi\mapsto J_V^*(s, \chi)$ \BBB is Fr\'echet differentiable in $H$ for every $s\in V'$,\BL
\item[(J2)] \BBB $\displaystyle\left\|\frac{\partial J_V^*(s, \chi)}{\partial \chi}\right\|_H\leq c_1\|\eta\|_H +c_2$, for every $\eta\in \partial_{V',V} J_V^*(s,\chi)$ and $(s,\chi)\in \EEEtris D(J_V^*)$, \BL 
\end{description}
where  $\frac{\partial(\cdot)}{\partial\chi}$ denotes the partial derivative with respect to $\chi$ (which will be denoted also by $\partial_\chi(\cdot)$ and by $(\cdot)_\chi$ in the paper). 

Moreover, we assume that 
\begin{description}
\item[(J3)] $ s\mapsto J_V^*(s,\chi)$ is proper, convex and lower semicontinuous from $V'$ to $[0,+\infty]$, for every $\chi\in H$,
\end{description}
so that the subdifferential  $\partial_{V',V} J_V(s,\chi)$ which maps $V'\times H$ into $2^V$ according to the definition \eqref{partialV}
turns out to be a maximal monotone operator acting from $V'$ to $2^V$, for every $\chi\in H$ (cf.~\cite{barbu}). \BL
\enhy
Note that, the assumption {\bf (J3)} follows from assumptions on $j$ and \eqref{JstarV} and  that possible examples of functions $j$ complying with our assumptions will be listed in the next Subsection~\ref{ex:jstar}. 
\beos\label{lowbou}
Observe that the assumptions on the positivity of the maps $j$ and $J_V^*$ could be weakened: we need indeed to have only a lower bound (possibly with a negative constant) for them in order to perform the first a-priori estimate \eqref{estI}. Moreover, let us nothe that the assumption {\bf (J2)} could be relaxed: we could indeed assume $c_1$, $c_2$ to be two continuous functions of $\chi$ bounded on ${\rm dom}(\widehat\beta)$.  However, we put ourselves in this setting to avoid further technicalities for the reader's convenience.
\eddos

\beos\label{domR} 
Note that, \BL under particular assumptions on the function $j$ (for example in case ${\rm dom}(j)=\RR$), we could also rewrite the functional $e$ in \eqref{inten} as (cf. Remark~\ref{rem:domR} for more details  and \cite{barbu}\BL) 
\begin{equation}\label{intenbis}
e(s,\chi,\nabla\chi)=\io\left(j^*(s, \chi)+\frac12|\nabla\chi|^2+W(\chi)\right)\dd x\,,
\end{equation}
where $j^*$ is the conjugate function of $j$ with respect to the variable $s$, i.e. 
\[
j^*(s,\chi)=\sup_{\theta\in \RR}\left(s\theta-j(\theta, \chi)\right),\quad \forall (s,\chi)\in \RR\times\RR\,. 
\]
\eddos 
We introduce as dissipative variables the time derivatives $s_t$ and $\chi_t$ \EEEbis (see, e.g. \cite{fremondlibro}, for a  definition of the pseudo-potential of dissipation {\sl \`a la Moreau})  and we include \BL dissipation in the model by choosing the following form for the pseudopotential of dissipation 
 depending on \BL the dissipative variables $s_t$ and $\chi_t$.  Note that \BL we suppose  the evolution to be rate dependent. The first possibility we  consider for \BL the dissipation functional  is 
\begin{equation}\label{phi}
\varphi:V'_0\times H\rightarrow\RR,\quad
\varphi(s_t,\chi_t)=\frac 1 2\duav{s_t, A^{-1}s_t}+\frac12\io|\chi_t|^2\dd x\,.
\end{equation} \BL
Note that \EEEbis we have a natural scalar product in $V_0'$ $((\cdot,\cdot))_{V_0'}$ defined as  $((s_t,s_t))_{V_0'}=\duav{s_t, A^{-1}s_t}=\|s_t\|_{V_0'}^2$, \BL  $\io|\chi_t|^2\dd x=(\chi_t,\chi_t) =\|\chi_t\|_H^2$. \BL
In the defintions of $e$ and $\varphi$ we have normalized all the physical constants to $1$ for simplicity and without any loss of generality. 
Another \BL  possibility consists in letting
\begin{equation}\label{phidue}
\varphi(s_t,\chi_t)=\frac 1 2\duav{s_t, A^{-1}s_t}+\frac 1 2\duav{\chi_t, A^{-1}\chi_t},
\end{equation} 
but we prefer not to exploit this case in the present contribution in order to distinguish between  the roles of the two variables: the thermal variable $s$ ({\sl conserved}) and the mechanical variable $\chi$ ({\sl non conserved}). Moreover, a rate-independent model could be introduced in place of the rate-dependent one we analyze here by suitably modifying the choice of the dissipation functional \eqref{phi} (cf., e.g.,~\cite{mt}). However, the analysis we are performing does not apply to this case, which would require ad hoc techniques and some suitable notion of weak solution.

\paragraph{The constitutive relations and the PDEs.}
\EEE
Now, according to the definition of $\mathcal{P}_{int}$  and of $e$ and $\varphi$ (cf.~\eqref{inten} and \eqref{phi}), \BL we let the thermal force $F$ be \BL 
\begin{equation}\label{defF}
 F=\partial_s e+\partial_{s_t}\varphi=\partial_{V',V} J_V^*(s,\chi)+A^{-1}(s_t)\,.
\end{equation}
Hence, for the evolution of $\chi$ 
we prescribe the following mechanical (micro) forces and stresses $B$ and $E$: \BL 
\begin{equation}
\label{vpchi2}
B=\partial_\chi e+\partial_{\chi_t}\varphi=\partial_\chi J_V^*(s,\chi)+\partial\widehat\beta(\chi)+\widehat\gamma'(\chi)+\chi_t, \quad E=\partial_{\nabla\chi} e=\nabla\chi\,.
\end{equation}
\BL

\EEE From \BBBnew (\ref{vps}--\ref{vpbouchi}) \EEE and the above constitutive \BBBnew relations  \BL we deduce the following \BBBnew PDE system  \BL for the evolution of $s$ and $\chi$:
\begin{align}\label{evs}
&s_t+\BBB A\eta=0 \,\hbox{in }V_0',\,\EEEtris\eta \in \partial_{V',V} J_V^*(s,\chi),\quad\hbox{a.e. in }(0,T)\BL\,,\\
\label{evchi}
&\chi_t-\Delta\chi+\BBBnew \xi\BL+\gamma(\chi)+\partial_\chi J_V^*(s,\chi)= 0, \,\EEEtris\xi\in \beta(\chi)\quad\hbox{a.e. in }\Omega\BBBnew\times(0,T)\BL\,\\
\label{evbouchi}
&\nabla\chi\cdot{\bf n}=0 \quad \hbox{ a.e. on }\Gamma\BBBnew \times(0,T)\BL\,,
\end{align}
where we denote by $\beta$ the subdiferential of $\widehat\beta$ ($\beta=\partial\widehat\beta$) and by $\gamma=\widehat\gamma'$.

Notice that system  \eqref{evs}--\eqref{evchi} can be rewritten in terms of the vector $u:=(s,\chi)$ in a more general framework, as the gradient-flow associated to the functional 
\begin{equation}
\label{defPhi}
\Phi(u)=\Phi(s,\chi)=J_V^*(s,\chi)+\io\left(\frac12|\nabla\chi|^2+W(\chi)\right)\dd x\,,
\end{equation}
as 
\begin{equation}\label{eqU}
\mathcal{N}(u_t)+\frac{\delta\Phi}{\delta u}\EEEtris\ni\BL0\quad \hbox{in }(0,T)\,,
\end{equation}
where $\mathcal{N}$ is the duality map between \EEEbis ${\cal H}:={V}_0'\times H$ \BL and $V_0\times H$ induced by the norm 
\begin{equation}\label{defN}
\|u\|_{\mathcal{H}}:=\duav{A^{-1}(s), s}+\io|\chi|^2\, \dd x, \quad \hbox{so that }  \mathcal{N}(s,\chi):=(A^{-1}(s), \chi)\,.
\end{equation}

\beos\label{rem:domR} Let us notice that in case ${\rm dom }(j)=\RR$, which is also equivalent to \EEE assume \BL
\[
\lim_{|r|\to+\infty} \frac{j^*(r)}{|r|}=+\infty\,,
\]
then we can prove that (cf.~\cite{BCGG}) the functional $e$ defined in \eqref{inten} can be rewritten as 
\[
e(s,\chi,\nabla\chi)=\io\left(j^*(s, \chi)+\frac12|\nabla\chi|^2+W(\chi)\right)\dd x\,,
\]
where $j^*$ is the conjugate function of $j$ with respect to the variable $s$, i.e. 
\[
j^*(s,\chi)=\sup_{\theta\in \RR}\left(s\theta-j(\theta, \chi)\right),\quad \forall (s,\chi)\in \RR\times\RR\,. 
\]
Moreover, the inclusion \eqref{evs} can be rewritten as the following gradient flow in $V_0'$:
\[
s_t+\partial_{V'}J_V^*(s,\chi)\ni 0\,,
\]
where $\partial_{V'} J_V^*$ is defined as the subdifferential of $J_V^*$ in $V'$ mapping $V'\times H$ into $2^{V'}$ as follows:
\begin{align}\label{partialVp}
&\xi\in \partial_{V'}J_V^*(s,\chi) \quad \hbox{ if and only if } \xi\in V', \\
\no
&\qquad (s, \chi)\in D(J_V^*),\hbox{ and } J_V^*(s,\chi)\leq (\!(s-w,\xi)\!)_*+J_V^*(w,\chi)\quad \forall (w,\chi)\in V'\times H\,,
\end{align}
where $(\!(\cdot,\cdot)\!)_*$ denotes the scalar product in $V'$. The reader can refer to \cite{barbu} and to \cite[Section~2]{BCGG} for the proofs of these results. Finally in this case we have $u\in \partial_{V',V} J_V(s, \chi)$ in $V$ iff $u\in \partial_s j^*(s, \chi)$ a.e. in $\Omega$, where $\partial_s$ denotes here the subdifferential of convex analysis with respect to the variable $s$ (cf., e.g., \cite{brezis}). 
\eddos
 
\subsection{Possible choices of $j^*$}
\label{ex:jstar}

In this section we show how to derive different types of phase-field models by our general system. \BL

\paragraph{The Caginalp model of phase transitions.} Choose $j^*(s,\chi)=\frac{s^2}{2}-s\chi+\frac{\chi^2}{2}$. Denote by $\theta:=\partial_s j^*=s-\chi$. Then, the Hyp.~\ref{HJ} is obviously satisfied and \BL the PDEs (\ref{evs}--\ref{evchi}) can be rewritten as 
\begin{align}\no
&\theta_t+\chi_t-\Delta\theta=0\,,\\
\no
&\chi_t-\Delta\chi+\beta(\chi)+\gamma(\chi)-\theta\ni0\,,
\end{align}
coupled with \BL Neumann homogeneous boundary conditions on $\theta$ and $\chi$\BL. 
This PDE system  can be easily identified with the ``standard'' phase field model of Caginalp type (cf.~\cite{caginalp}), letting $\theta$ be the relative temperature of the system and $\chi$ the local proportion of one of the two phases of the substance undergoing phase transitions\BL .

\paragraph{The entropy model for phase transitions.} Choosing $j^*(s,\chi)=j^*(s-\lambda(\chi))=\exp(s-\lambda(\chi))$, we have that Hyp.~\ref{HJ} is satisfied in case $\lambda$ is a Lipschitz continuous function on the domain of $\beta$. Then, defining
$\theta:=\partial_sj^*=\exp(s-\lambda(\chi))$, we get  \BL $s=\log\theta+\lambda(\chi)$ and 
the PDEs (\ref{evs}--\ref{evchi}) can be rewritted as 
\begin{align}\no
&(\log\theta+\lambda(\chi))_t-\Delta\theta=0\,,\\
\no
&\chi_t-\Delta\chi+\beta(\chi)+\gamma(\chi)-\lambda'(\chi)\theta\ni0\,,
\end{align}
again with Neumann homogeneous boundary conditions for both $\theta$ and $\chi$. \BL 
This system can be easily identified with the ``entropy'' phase field model introduced in \cite{bcf} and \cite{BFI}. Here $\theta$ respresents the absolute temperature of the system \BL which is forced to be positive, by the presence of the logarithmic nonlinearity in the $\theta$-equation\BL. Let us notice that in this case the assumption $D(j^*)=\RR$ is not verified, hence we are not entitled to use the function $j^*$ instead of the operator $J_V^*$ in $e$ (cf. Remark~\ref{domR}), so, the choice we made here is only formal. For a rigorous analysis of this case the reader can refer to \cite{BFR}.

\paragraph{The Penrose-Fife model for phase transitions.}  We choose \EEE $j^*(s,\chi)=-\log(s-\chi)$,  for $s>\chi$ and define \EEE $\theta:=\partial_sj^*(s,\chi)=-\frac{1}{s-\chi}$. Then, \EEE we observe that we can formally get the Penrose-Fife mode. \EEE Indeed, it results that $\partial_\chi  j^*(s,\chi)=\frac 1 {s-\chi}=\frac 1 \theta$. Thus, we can rewrite (\ref{evs}--\ref{evchi}) as 
\begin{align}\no
&(\theta\pm\chi)_t\mp\Delta\left(-\frac 1\theta\right)=0\,,\\
\no
&\chi_t-\Delta\chi+\beta(\chi)+\gamma(\chi)+\frac{1}{\theta}\ni0\,,
\end{align}\BL
coupled with \BL Neumann homogeneous boundary conditions on \EEEbis  $-\frac1\theta$ \BL and $\chi$\BL.  
This system can be easily identified with the Penrose-Fife model of phase transitions introduced in \cite{pf}. 

\beos\label{case2}
Let us notice that in case we choose as \BL pseudopotential of dissipation \EEEbis the functional \eqref{phidue}   in \eqref{defF} and \eqref{vpchi2}, \BL the first equation is the same as \eqref{evs}, while the equation for $\chi$ results as:\BL
\begin{equation}
\label{evchiab}
A^{-1}\BL\chi_t-\Delta\chi+\beta(\chi)+\gamma(\chi)+\partial_\chi j^*(s,\chi)\ni0\,,
\end{equation}
and thus
\begin{align}\label{evsch}
&s_t+A\left(\partial_s j^*(s,\chi)\right)\ni0 \quad\hbox{in }V_0'\BL \,,\\
\label{evchich}
&\chi_t+A w=0, \quad w\in-\Delta\chi+\beta(\chi)+\gamma(\chi)+\partial_\chi j^*(s,\chi)\,,
\end{align}
coupled with Neumann homogeneous boundary conditions for $\partial_s j^*(s,\chi)$, $\chi$ and $w$. In this case the evolution of $\chi$ is ruled by the well-known fouth order Cahn-Hilliard equation modelling phase separation phenomena (cf., e.g., \cite{CH}). However, as we already mentioned, we prefer not to deal with this case here. 
\eddos

%%%%%%%%%%%%%%%%%%%%%%%%%%%%%%%%%%%%%%%%%%%%% 
\section{Main results} 
\label{sec:main}

In this section we state the main results of this paper, the first one  (Thm.~\ref{th1}) concerns the existence of global in time solutions for system \eqref{evs}--\eqref{evchi} coupled with the boundary condition \eqref{vpbouchi} and the initial conditions
\begin{align}\label{inis}
&s(0)=s_0 \quad \hbox{in }\EEE D(J_V^*)\BL\,,\\
\label{inichi}
&\chi(0)=\chi_0\quad \hbox{a.e. in }\Omega\,,
\end{align}
while the second one (Thm.~\ref{th2}) regards uniqueness of solutions under more restrictive assumptions on the nonlinearities involved. Let us start with the first result. 
\BL

\bete\label{th1} Assume Hypotheses~\ref{HJ} and \ref{HW} and take $s_0\in D(J_V^*)$, $\chi_0\in V\cap {\rm dom}(\widehat\beta)$. Then, for every $T>0$ there exists at least one solution $(s,\chi)$ to \eqref{evs}--\eqref{evbouchi} \BBBnew and \eqref{inis}--\eqref{inichi} \BL satisfying the regularity properties:
\begin{align}\label{regos}
&s\in H^1(0,T;V_0')\,,\\
\label{regochi}
&\chi\in H^1(0,T;H)\cap L^\infty(0,T;V)\EEE\cap L^2(0,T;H^2(\Omega))\BL\,.
\end{align}
\ente

\prova 
In order to prove Theorem~\ref{th1}, we first approximate system (\ref{evs}--\ref{evchi}) with a regularized problem depending on a positive small parameter $\e$ and then \EEE we pass to the limit by (weak-strong) compactness arguments and semicontinuity results  based on \BL sufficient a-priori estimates -- independent of $\e$ -- \EEE we are going to prove on the approximating solutions. \BL

\paragraph{The approximated problem.}  \EEE Let us fix $\varepsilon>0$. Then, \BL  for every $T>0$ and $(s_{0, \e}\BL,\chi_0)\in (D(J_V^*)\cap H)\times (V\BBBnew \cap {\rm dom}(\widehat\beta)\BL) $,  \EEE we aim to find \BL a solution $(s_\e,\chi_\e)\in H^1(0,T; V_0'\times H)$  \EEE to the following \BBBnew differential inclusions:  \BL
\begin{align}\label{evse}
&\dt s_\e+\BBB  A\left(\eta_\e+\e s_\e\right)= 0\,\EEEtris\hbox{ in }V_0',\, \eta_\e\in \partial_{V',V} J_V^*(s_\e,\chi_\e),\quad \hbox{ for a.e. } t\in (0,T)\,,\BL\\
\label{evchie}
&\dt\chi_\e-\Delta\chi_\e+\BBB \xi_\e+\gamma(\chi_\e)+\partial_\chi J_V^*(s_\e,\chi_\e)=0\quad\xi_\e\in \beta(\chi_\e),\quad \hbox{a.e. in }\Omega\times (0,T)\BL\,,
\end{align}
coupled with the boundary and initial conditions \eqref{evbouchi} and (\ref{inis}--\ref{inichi}), with $s_{0,\e}$ in place of $s_0$. \EEE In particular, we
assume that 
\begin{equation}\label{inizse}
\BBBnew s_{0, \e}\in D(J_V^*)\cap H, \quad s_{0,\e}\to s_0\quad\hbox{ in }V_0'\quad\hbox{as $\e\searrow0$}.
\end{equation}
\BL

\EEE We first observe that we can recover \eqref{evse} and and \eqref{evchie}, by approximating the energy functional  \eqref{defPhi} as follows: 
\[
\Phi_\e(s,\chi):=\Psi_\e(s,\chi)+\io\left(\frac12|\nabla\chi|^2+W(\chi)\right)\dd x\,, \,\, \Psi_\e(s,\chi):=J_V^*(s,\chi)+\frac{\e}{2}\io|s|^2\dd x.
\]
\EEE Actually, note that now $\Phi_\e$ is defined in $(V_0'\cap H)\times H$.  Hence, we can construct its subdifferential in the duality between $V'_0$ and $V_0$, and 
rewrite the \EEE equation \eqref{evse} as
\begin{align}\label{evsePhi}
&\dt s_\e+\BBB A\zeta_\e=0 \,\BL\hbox{ in }V_0',\,\EEEtris\zeta_\e\in \partial_{V',V} \Phi_\e(s_\e,\chi_\e)\quad\quad\hbox{ for a.e. } t\in (0,T)\,.\BL
\end{align}
\EEE Now, our aim is to 
prove the existence of solutions \BBBnew of \eqref{evse}--\eqref{evchie}, \eqref{inis}--\eqref{inichi} with $s_{0,\e}$ istead of $s_0$, and \eqref{evbouchi} by a time-discrete approximation, as follows (cf.~also \cite{RS} for a similar procedure). Here we drop the index $\e$ in order to simplify the notation. 
Let us fix a time step $\tau=T/N$, $N\in \NN$ and introduce a uniform partition
\[
P_\tau:=\{t_0=0,\, t_1=\tau,\,\dots, \, t_n=n\tau, \, \dots, \, t_N=T\}
\]
of the interval $(0,T)$. Then, we need to find a discrete approximation $s^n\sim s(t_n)$,  $\chi^n\sim\chi(t_n)$ by solving the implicit Euler scheme (cf.~also \eqref{eqU}):
\begin{equation}\label{euler}
\mathcal{N}\left(\frac{U^n-U^{n-1}}{\tau}\right)+\EEEtris\zeta^n=0, \BL\quad n=1, \, \dots, \, N\,; U^0:=u_0\,,
\end{equation}
where \EEEtris $\zeta^n\in\frac{\delta\Phi_\e}{\delta u}(U^n)$ and \BL we have defined $U^n=(s^n,\chi^n)$, $u_0=(s_{0,\e},\chi_0)$. \EEE Using the functional space, \EEEbis we have already introduced to define the operator $\mathcal{N}$, \EEE ${\cal H}=(V_0'\cap H)\times H$, we \BL notice that \eqref{euler} is the Euler equation for the variational problem 
\begin{equation}
\label{var}
\begin{cases}
&\hbox{find } U^n\in \mathcal{H}\quad\hbox{minimizing}\\
&F_\e(\tau, U^{n-1}; U):=\frac{1}{2\tau}\left\| U-U^{n-1}\right\|_{\mathcal{H}}^2+\Phi_\e(U), \quad U\in \mathcal{H}\,.
\end{cases}
\end{equation}
It is not difficult to see that this minimization problem is solvable due to the lower-semicontinuity and coercivity properties of $\Phi_\e$ (cf., e.g., \BBB\cite{RS, RSESAIM} \BL and references therein for a similar variational approach to find a discrete solution\BL). 

Then, we can construct the piecewise constant interpolants $\bar{U}_\tau(t):=U^n$ if $t\in ((n-1)\tau, \, n\tau]$.  In particular, we get that  and we  recover the solution $U:=(s,\chi)$ of (\ref{evse}--\ref{evchie})   as the limit of $\bar{U}_\tau$ as $\tau\searrow0$.This can be done, by using suitable a priori estimates (independent of $\tau$ and then passing to the limit by compactness and semicontinuity arguments. We do not enter the details of the proof, as it is very similar to the estimates and passage to the limit procedure we are going to detail in the next sections to pass to the limit as $\e\searrow0$.   Note that  in this case some technicalities are avoided due to the more regular setting for the variable $s_\e$ \EEEbis (recall  the strict positivity of $\e$).
\EEE
Thus, we \BBBnew can easily \EEE prove the following existence result. 
\bete
Under the same assumptions of Theorem \ref{th1}, letting $\e>0$ be fixed and \eqref{inizse} holds, then there exists a solution to \eqref{evse}-\eqref{evchie} with $s_\e(0)=s_{0,\e}$ and $\chi_\e(0)=\chi_0$, with the following regularity
\begin{align}\label{regs}
&s_\e\in H^1(0,T;V_0')\cap L^\infty(0,T;V_0),\\\label{regchie}
&\chi_\e\in H^1(0,T;H)\cap L^\infty(0,T;V)\cap L^2(0,T;H^2(\Omega)).
\end{align}
\ente
\BL

%%%%%%%%%%%%%%%%%%%%%%%%%%%%%%%%%%%%%%%%%%%%%%%%%%%%%%%%%%%%%%%%%%%%%%%%%%%%%%%%%%%%%%%%%%%%%%%%%%%%%%
\paragraph{A priori estimates (uniform in $\e$).}
\label{sec:estimate}

Let us consider the system (\ref{evse}--\ref{evchie}), where \EEE for the sake of simplifying notation \BL we neglect the index $\e$ \EEE for solutions and involved functions.\BL  We now perform the a priori estimates independent of $\e$, so, we use here the same symbol $c$ for positive constants, possibily different from line to line, depending on the problem data, but independent of $\e$. 

\BBB
In order to perform \EEEtris the \BBB first estimate we need to prove here a small variant of the chain rule formula stated, e.g., in 
 \cite[Prop.~4.2]{ckrs}. 
\bepr\label{chain}
Let $G: V'\times H \to [0,+\infty]$ be a map such that 
\begin{description}
\item $v\mapsto G(u,v)$ is Fr\'echet differentiable for every $u\in V'$,
\item $u\mapsto G(u,v)$ is a proper convex lower semicontinuous mapping for every $v\in H$,
\end{description}
and let $u\in H^1(0,T;V')\cap L^2(0,T;V_0)$, $v\in H^1(0,T;H)\cap
L^2(0,T; V)$, and
$\delta(t) \in \partial_{V',V}G(u(t), v(t))$ for a.e. $t\in (0,T)$, where the subdifferential $\partial_{V',V}$ is defined as in \eqref{partialV}.
Then the function $g=
G(u(\cdot), v(\cdot))$ is absolutely continuous in $[0,T]$ and
$g'(t) = \duav{u'(t),\delta(t)}+(v'(t), \partial_v G(u(t),v(t)))$ for a.e. $t\in (0,T)$.
\empr

\prova
Here we follow the lines of the proof of \cite[Prop.~4.2]{ckrs}. 
Let $w \in W^{1,\infty}(0,T)$ be a non-negative
function with compact support in $(0,T)$. We choose $h>0$ such
that supp$ {(w) \subset [h, T-h]}$. For a.e. $t \in\EEEtris(0,T)$, by definition of sub-differentials we can infer that   \BBB  
\begin{align}\no
& {\duav{u(t) - u(t-h), \delta(t-h)}}+(v(t)-v(t-h), \partial_v G(u(t-h), v(t-h))\\
\no
& \leq g(t) - g(t-h) \leq \duav{ u(t) - u(t-h), \delta(t)}+(v(t)-v(t-h), \partial_v G(u(t), v(t)).
\end{align}
\EEEtris Indeed, observe that $(\delta,\partial_v G)$ belongs to $\partial G$, i.e. to the sub-differential of the function $G:V'\times H\rightarrow[0,+\infty]$ defined w.r.t.~the variable $(u,v)$.  
Observe that we can extend $w$ outside of $(0,T)$ with the $0$ value.
Hence, multiplying by $w(t)$, integrating with respect to $t$, and
letting $h\searrow 0$, we obtain
\begin{align}
\no
\frac 1h \int_h^T \duav{ u(t) - u(t-h) , \delta(t-h) }w(t)dt
& = \frac 1h \int_0^{T-h} \! \duav{ u(t+h) - u(t) , \delta(t) }w(t+h)dt\\
\no & \to \ \int_0^T \duav{u'(t),\delta(t)}w(t)dt,\\
\no \frac 1h \int_h^T (g(t) - g(t-h))w(t)dt&=
\frac 1h \int_0^T \delta(t)(w(t)- w(t+h))dt\\
\no &\to \ -\int_0^T \delta(t)w'(t)dt,\\
\no \frac 1h \int_h^T \duav{ u(t) - u(t-h) , \delta(t) }w(t)dt
& \to\ \int_0^T \duav{u'(t), \delta(t) }w(t)dt.
\end{align}
Moreover, 
\begin{align}
\no
&\frac 1h \int_h^T(v(t)-v(t-h), \partial_v G(u(t-h), v(t-h))w(t)dt\\
\no
&= \frac 1h \int_0^{T-h} (v(t+h) - v(t) , \partial_v G(u(t), v(t))w(t+h)dt\to \ \int_0^T (v'(t), \partial_v G(u(t), v(t))w(t)dt,\\
&\no \frac 1h \int_h^T (v(t)-v(t-h), \partial_v G(u(t), v(t))w(t)dt
\to\ \int_0^T(v'(t), \partial_v G(u(t), v(t))w(t)dt.
\end{align}
Therefore, we conclude that
$$
-\int_0^T g(t)w'(t)dt \ = \ \int_0^T \left(\duav{u'(t), \delta(t)}+(v'(t), \partial_v G(u(t), v(t)))\right)w(t)dt
$$
for all non-negative Lipschitz continuous test functions $w$
with compact support. Since both the positive and the negative
part of a Lipschitz continuous function are Lipschitz
continuous, we obtain the assertion.
\qed

\BL

\noindent{\sl First a priori estimate.}\quad We test \eqref{evse} by $A^{-1}s_t$ getting
\begin{equation}
\langle s_t,A^{-1}s_t\rangle=\|s_t\|^2_{V_0'},
\end{equation}
and in addition, \BBBnew using the definition of $J_V^*$ and a variant of the chain rule formula stated in Proposition~\ref{chain} with the choices $G=J_V^*$, $u=s$, $v=\chi$, we get 
\begin{align}
&\langle A\BBB \eta,A^{-1}s_t\rangle=\duav{s_t,\eta}=\frac d{dt}J_V^*(s(t),\chi(t))-(\chi_t,\partial_\chi J_V^*(s,\chi))\, .
\end{align}
\BL
Testing \eqref{evchie} by $\chi_t$, we get 
\begin{equation}
\|\chi_t\|^2_H+\frac 12\frac d{dt}\|\nabla\chi\|^2_H+\frac d{dt}\int_\Omega W(\chi)+\BBB (\chi_t,\partial_\chi J_V^*(s,\chi))\BL=0\,.
\end{equation}
Adding the resulting equations and integrating over $(0,t)$, $t\in (0,T)$,  \BBB and using  the definition of $\Psi_\e$, \BL we obtain 
\begin{equation}\label{estI}
\int_0^t\left(\|s_t\|^2_{V_0'}+\|\chi_t\|^2_H\right)\dd \tau+\Psi_\e(s(t),\chi(t))\BL+\int_\Omega W(\chi(t))+\|\nabla\chi(t)\|^2\leq c,
\end{equation}
where here $c$ depends in particular on the initial data. 
Adding to both sides in \eqref{estI} 
$$\|\chi(t)\|_H^2=\|\chi_0\|_H^2+2\int_0^t(\chi(\tau), \chi_t(\tau))\, \dd \tau\,,$$ 
and using H\"older and Young inequalities together with  Hyp.~\ref{HW} and \ref{HJ} and  a standard Gronwall lemma\BL, we obtain
\begin{align}\label{stI}
&\|s_t\|^2_{L^2(0,T;V_0')}+\e\|s\|^2_{L^\infty(0,T;H)}\BL\leq c,\\\label{stII}
&\|\chi\|_{H^1(0,T;H)\cap L^\infty(0,T;V)}\leq c,\\\label{stIII}
&\|\widehat\beta(\chi)\|_{L^\infty(0,T;L^1(\Omega))}\leq c\,.
\end{align}

\noindent{\sl Second a priori estimate.}\quad We proceed by a comparison in \BBB \eqref{evsePhi}. \BL Due to \eqref{stI}$_1$, \BL we have that \BBB $A\zeta_\e$ \BL is bounded in $L^2(0,T;V_0')$, and thus 
\begin{equation}\label{stIV}
\BBB \|\zeta_\e\|_{L^2(0,T;V_0)}\leq c\,. \BL
\end{equation}
Hence, using \eqref{stI}$_2$, we get 
\begin{equation}\label{stV}
\|\BBB\eta_\e\BL\|_{L^2(0,T;H)}\leq c\,. 
\end{equation}
\EEE Indeed, due to the definition of $\Psi_\e$ we can infer that \BBB $\zeta_\e=\eta_\e+\partial \psi_\e(s)$, \BL where we have used the notation $\psi_\e(s)=\frac {\e}2\int_\Omega s^2{\rm d}x $ and the fact that $\psi_\e$  has sa domain the whole real line and thus its subdifferential in the duality $V',V$ corresponds to the standard subdifferential of the convex analysis (cf.~\cite{brezis}).\BL

\noindent
{\sl Third a priori estimate.}\quad Using now Hyp.~\ref{HJ} together with \eqref{stV}, we get 
\begin{equation}\label{stVI}
\|\partial_\chi J_V^*(s,\chi)\|_{L^2(0,T;H)}\leq c_1\|\BBB \eta_\e \BL\|_{L^2(0,T;H)}+c_2\leq c\,. 
\end{equation}
Moreover, by comparison in \eqref{evchie} and by standard monotonicity and regularity results, we get 
\begin{equation}\label{stVII}
\|\BBB\xi\BL\|_{L^2(0,T;H)}+\|\chi\|_{L^2(0,T;H^2(\Omega))}\leq c\,.
\end{equation}
\BL

\paragraph{Passage to the limit as $\e\searrow 0$. }
\label{sec:limit}

Now, we aim to pass to the limit in \eqref{evse}--\eqref{evchie} as $\e\searrow 0$, revovering finally a solution to (\ref{evs}--\ref{evchi}). By virtue of (\ref{stI}--\ref{stVII}) and by compactness results \BL we get (at least for some subsequences of $\e\searrow 0$):
\begin{align}\label{convI}
&\sn\weakstar s\quad\hbox{in }H^1(0,T;V_0'),\\\label{convII}
&\chin\weakstar\chi\quad\hbox{in }H^1(0,T;H)\cap L^\infty(0,T;V)\cap L^2(0,T;H^2(\Omega))\BL\,,\\\label{convIII}
&\BBB \zeta_\e\BL\weak \BBB \zeta\BL \quad\hbox{in }L^2(0,T;V_0)\,,\\\label{convIIIbis}
&\BBB \eta_\e\BL \weak \BBB \eta \BL \quad \hbox{in }L^2(0,T;H)\,,\BL\\\label{convIV}
&\partial_\chi J_V^*(\sn,\chin)\weak z \quad\hbox{in }L^2(0,T;H)\,,\\\label{convV}
&\e^{1/2}\BL \sn\weakstar 0 \quad \hbox{in } L^\infty(0,T;H)\,.
\end{align}
Notice that, by the definition of $\Psi_\e$ and by \eqref{convIII}, \eqref{convIIIbis}, and \eqref{convIV}, we immediately deduce that \BBB $\zeta=\eta$ \BL a.e..
Moreover, \BL by strong compactness, from \eqref{convII} we can also deduce (at least) (cf.~\cite{si})\BL
\begin{equation}\label{convVI}
\chin\rightarrow\chi\quad\hbox{in }C^0([0,T];H)\cap L^2(0,T;V)\BL\,.
\end{equation}
We aim to  identify \BBB $\zeta\in\partial_{V',V}J_V^*(s,\chi)$ \BL (see \eqref{convIII}). By definition of sub differential this corresponds to prove that
\begin{equation}\label{ideta}
\int_0^T\langle v-s,\BBB\zeta\BL\rangle\,\dd\tau\leq \int_0^T\left(J_V^*(v,\chi)-J_V^*(s,\chi)\right)\dd\tau\quad\forall v\in V_0'\quad\hbox{and }\chi\in H\,.
\end{equation}
\BL
Note that, if we test the equation \eqref{evse} by $A^{-1}\sn$ \BL and integrate in time, by weak lower semicontinuity of 
norm we have (for $\BBB\zeta_\e\in \partial_{V',V}\Psi_\e(\sn,\chin)$\BL)
\begin{equation}
\int_0^T\limsup_{\e\BL\searrow 0}\,\langle \sn,\BBB \zeta_\e\BL \rangle\dd\tau\leq \int_0^T\langle s,\BBB\zeta\BL\rangle\dd\tau\,.
\end{equation}
Hence, by employing \eqref{convIII}, \eqref{convIIIbis}, and \eqref{convV}, \EEE and using the fact that \BBB $\zeta_\e$ \BL  belongs to $\partial_{V',V}\Psi_\e(s_\e,\chi_\e)$,\BL  we get, for all $v\in V_0'$ and $\chi_\e\in H$, 
\begin{align}\label{liminfI}
&\int_0^T\langle v-s,\BBB\zeta\BL\rangle\dd\tau\leq \int_0^T\langle v,\BBB\zeta\BL\rangle\dd\tau-\limsup_{\e\searrow0}\int_0^T\langle \sn,\BBB \zeta_\e\BL\rangle\dd\tau\\\no
&=\int_0^T\langle v,\BBB\zeta\BL\rangle\dd\tau+\liminf_{\e\searrow0}\int_0^T\left(-\langle \sn,\BBB\zeta_\e\BL\rangle\right)\dd\tau\\\no
&\leq\liminf_{\e\searrow0}\int_0^T\left(\langle v,\BBB\zeta_\e\BL\rangle-\langle \sn,\BBB\zeta_\e\BL\rangle\right)\dd\tau\\\no
&\leq\limsup_{\e\searrow 0}\int_0^T\left(\Psi_\e(v,\chin)-\Psi_\e(\sn,\chin)\right)\dd\tau\\
\no
&=-\liminf_{\e\searrow0}\int_0^T\left(-J_V^*(v,\chin)+J_V^*(\sn,\chin)-\frac{\e}{2}\|v\|_H^2+\frac{\e}{2}\|s_\e\|_H^2\right)\dd\tau\\
\no
&\leq-\liminf_{\e\searrow0}\int_0^T\left(-J_V^*(v,\chin)+J_V^*(\sn,\chin)\right)\dd\tau \,.
\end{align}
\BL
Using Hyp.~\ref{HJ} and  \eqref{convVI}, we can deduce
\begin{equation}
\lim_{\e\searrow 0}J_V^*(v,\chin)=J_V^*(v,\chi)\,.
\end{equation}
Hence, let us observe that
\begin{equation}
J_V^*(\sn,\chin)-J_V^*(s,\chi)=J_V^*(s_\e,\chi_\e)-J_V^*(s_\e,\chi)+J_V^*(s_\e,\chi)-J_V^*(s,\chi)\BL,
\end{equation}\EEE
where
$$
I_1:=J_V^*(s_\e,\chi_\e)-J_V^*(s_\e,\chi),
$$
and 
$$
I_2:=J_V^*(s_\e,\chi)-J_V^*(s,\chi).
$$
Hence,
\BL 
by \eqref{stVI} and \eqref{convVI}, we have \BL
\begin{align}
&\int_0^T|I_1|\dd\tau\leq\int_0^T\|\partial_\chi J_V^*(\sn,\chin)\|_{H} \|\chin-\chi\|_{H}\dd\tau\rightarrow0\quad\hbox{as }\e\searrow 0\,.
\end{align}
The second integral
$$
\int_0^TI_2\,\dd\tau=\int_0^T\left(J_V^*(\sn,\chi)-J_V^*(s,\chi)\right)\dd\tau
$$
\BL
is treated by using the fact that $J_V^*$ is lower semicontinuous in $\sn$ with respect to weak convergence as it is convex (for $\chi$ fixed) and so 
\[
\liminf_{\e\searrow 0}\int_0^T I_2\dd\tau \geq 0
\]
Hence, coming back to \eqref{liminfI}, and using \eqref{HJ}, we have, for all $v\in V_0'$ and $\chi_\e\in H$,
\begin{align}\
\no
&-\int_0^T\liminf_{\e\searrow 0}\left(-J_V^*(v,\chi_\e)+J_V^*(s_\e,\chi_\e)\right)\dd\tau\\
\no
&\leq -\int_0^T\liminf_{\e\searrow0}(-J_V^*(v,\chi_\e))\dd\tau-\int_0^T\liminf_{\e\searrow0}\left(J_V^*(s_\e,\chi_\e)-J_V^*(s,\chi)+J_V^*(s,\chi)\right)\dd\tau\\
\no
&\leq \int_0^TJ_V^*(v,\chi)-\int_0^T\liminf_{\e\searrow0} I_1-\int_0^T\liminf_{\e\searrow0}I_2-\int_0^T J_V^*(s,\chi)\dd\tau\\
\no
&\EEE\leq \int_0^T\left(J_V^*(v,\chi)-J_V^*(s,\chi)\right)\dd\tau\,,
\end{align}
and this concludes the proof of \eqref{ideta}. Finally, using the previous convergences, where we have now identified $\BBB\zeta=\eta\in\partial_{V',V}J_V^*(s,\chi)\BL$ in $L^2(0,T;H)$, we can pass to the limit in the approximated system (\ref{evse}--\ref{evchie}) as well as in the corresponding boundary and initial conditions for $\e\searrow0$. This concludes the proof of Theorem~\ref{th1}.\qed
\bigskip

We conclude now with the last result of our paper concerning uniqueness of solutions for problem (\ref{evs}--\ref{evbouchi}), (\ref{inis}--\ref{inichi}). 
\bete\label{th2}
Assume Hypotheses~\ref{HJ} and \ref{HW} and take $s_{0}\in D(J_V^*)$, $\chi_{0}\in V\cap {\rm dom}(\widehat\beta)$ and suppose moreover that 
\begin{equation}\label{ipouni}
\hbox{the maps } \chi\mapsto \partial_{V,V'} J_V^*(s,\chi)\quad \hbox{and } \chi\mapsto \partial_\chi J_V^*(s,\chi)\hbox{ are Lipschitz continuous}
\end{equation}
from $H$ to $V_0$ and from $H$ to $H$, respectively, for every $s\in V'$. Then the solution $(s,\chi)$ of (\ref{evs}--\ref{evbouchi}), (\ref{inis}--\ref{inichi}) is uniquely determined and the following continuous dependence estimate holds true:
\begin{align}\label{CD}
\|(s_1-s_2)(t)\|_{V'}^2+\|(\chi_1-\chi_2)(t)\|_H^2+\int_0^t\|\nabla(\chi_1-\chi_2)\|_H^2\dd\tau\leq &C\left(\|(s_1-s_2)(0)\|_{V'}^2\right.\\
\no
&\left.+\|(\chi_1-\chi_2)(0)\|_H^2\right)\,.
\end{align}
\ente
\prova Let us take the difference of the two equations \eqref{evs} and the two relations \eqref{evchi} \EEE corresponding \BL to two different solutions $(s_i,\chi_i)$, $i=1,\,2$ and test them by $A^{-1}(s_1-s_2)$ and $(\chi_1-\chi_2)$, respectively. Integrating over $(0,t)$, for $t\in [0,T]$, using Hyp.~\ref{HJ} and Hyp.~\ref{HW}, we get 
\begin{align}\no
&\|(s_1-s_2)(t)\|_{V'}^2+\|(\chi_1-\chi_2)(t)\|_H^2+2\int_0^t\|\nabla(\chi_1-\chi_2)\|_H^2\dd\tau\leq 
\|(s_1-s_2)(0)\|_{V'}^2\\
\no
&\qquad+\|(\chi_1-\chi_2)(0)\|_H^2-2\int_0^t(\gamma(\chi_1)-\gamma(\chi_2), \chi_1-\chi_2)\dd\tau \\
\no
&\qquad-2\int_0^t\duav{s_1-s_2,\partial_{V,V'}J_V^*(s_2,\chi_1)-\partial_{V,V'}J_V^*(s_2,\chi_2)}\dd\tau\\
\no
&\qquad-2 \int_0^t(\chi_1-\chi_2,\partial_\chi J_V^*(s_2,\chi_1)-\partial_\chi J_V^*(s_2,\chi_2))\dd\tau\,.
\end{align}
Using then the Lipschitz continuity of $\gamma$ and assumption \eqref{ipouni}, together with Gronwall lemma, we obtain exactly \eqref{CD}. 
\qed

\BL

%%%%%%%%%%%%%%%%%%%%%%%%%%%%%%%%%%%%%%%%%%%%%%%%%%%%%%%%%%%%%%%%%%%%%%%%

\BBB

\section*{Acknowledgements}
The financial support of the FP7-IDEAS-ERC-StG \#256872
(EntroPhase) is gratefully acknowledged by the authors. The present paper 
also benefits from the support of the MIUR-PRIN Grant 2010A2TFX2 ``Calculus of Variations'' for EB, the GNAMPA (Gruppo Nazionale per l'Analisi Matematica, la Probabilit\`a e le loro Applicazioni) of INdAM (Istituto Nazionale di Alta Matematica), and the IMATI -- C.N.R. Pavia for EB and ER.
The authors would also like to thank Riccarda Rossi for the fruitful discussions on the topic. 
\BL
%%%%%%%%%%%%%%%%%%%%%%%%%%%%%%%%%

%\fer{La biblio \`e da ampliare  dopo aver scritto l'Introduzione}

\end{document}